\documentclass[10pt]{article}
\usepackage[francais]{babel}
\usepackage[ansinew]{inputenc}
\usepackage{amsmath}

\title{Cohomologie $L^p$ en degré 1 des espaces homog{\`e}nes
\footnote{Mots cl{\'e} : Cohomologie $L^{p}$, courbure n{\'e}gative, espace
homog{\`e}ne, espace de Besov. Keywords : $L^{p}$-cohomology, negative curvature, homogeneous space, Besov space.
Mathematics Subject Classification :
43A15, 
43A80, 
46E35, 
53C20, 
53C30, 
58A14. 
}}
\author{Pierre Pansu$^{1,2}$\footnote{$^{1}$ Univ Paris-Sud, Laboratoire de Mathématiques d'Orsay, Orsay, F-91405 ;\hfill\eject\indent\hskip7pt $^{2}$ CNRS, Orsay, F-91405.}
}

\newtheorem{theo}{Th{\'e}or{\`e}me}
\newtheorem{lemme}{Lemme}
\newtheorem{prop}[lemme]{Proposition}
\newtheorem{cor}[lemme]{Corollaire}
\newtheorem{defi}[lemme]{D{\'e}finition}

\def\preuve{\par\medskip\noindent {\bf Preuve.}{\hskip1em}}

\def\qed{~q.e.d.}

\def\R{{\bf R}}
\def\C{{\bf C}}

\def\ltimes#1{\times_{\alpha}}
\def\n#1{{\parallel #1 \parallel}}
\def\l{\lambda}

\def\N{{\cal N}}

\def\G{{\cal G}}
\def\A{{\cal A}}

\def\dim{{\rm dim}\,}
\def\ker{{\rm ker}\,}

\def\tr{{\rm tr}\,}

\def\con{\hbox{const.}}
\def\d{\displaystyle}
\def\inv#1{\frac{1}{#1}}

\def\hp#1{H^{#1,p}}

\def\tp#1{T^{#1,p}}

\def\rp#1{R^{#1,p}}

\def\lp{_{L^{p}}}

\def\timez{\times\hskip-2pt\vrule height 1ex width .055ex depth -.0ex}

\begin{document}
\maketitle
\begin{quote}
{\small
RESUME. On calcule la cohomologie $L^p$ en degré 1 des espaces homog{\`e}nes riemanniens. Notamment, on montre que la cohomologie r{\'e}duite est non nulle si et seulement si l'espace est quasiisom{\'e}trique {\`a} un espace homog{\`e}ne {\`a} courbure sectionnelle strictement n{\'e}gative.

ABSTRACT. The $L^p$-cohomology in degree 1 of Riemannian homogeneous spaces is computed. It turns out that reduced cohomology does not vanish exactly for spaces quasiisometric to negatively curved homogeneous spaces.}
\end{quote}

\section{Introduction}

L'objet de ce papier est de calculer dans des exemples, en employant des outils d'analyse, un invariant d'inspiration topologique, utile en théorie géométrique des groupes, la cohomologie $L^p$ en degré 1.

\subsection{Cohomologie $L^p$}

\begin{defi}
\label{hp1}
Soit $p\geq 1$. La cohomologie $L^p$ en degré 1 d'une variété riemannienne $M$, notée $\hp1(M)$, est le quotient de l'espace des 1-formes fermées $L^p$ par le sous-espace des différentielles des fonctions $L^p$ à gradient $L^p$. Autrement dit, $\hp1(M)=0$ si et seulement si toute fonction sur $M$ dont le gradient est $L^p$ est elle-même $L^p$ à une constante additive près. 
Les classes des 1-formes qui sont limites $L^p$ de différentielles de fonctions $L^p$ forment un sous-espace de $\hp1(M)$ appelé \emph{torsion} et noté $
\tp1(M)$. Le quotient $\hp1(M)/\tp1(M)=\rp1(M)$ (plus grand quotient séparé de $\hp1(M)$) s'appelle la \emph{cohomologie $L^p$ réduite}.
\end{defi}

Lorsque $M$ est compacte, $\hp1(M)$ est simplement le premier espace de cohomologie de de Rham de $M$. Dans cet article, on calcule $\hp1(M)$ lorsque $M$ est homogène non compacte. Voici deux cas extrêmes, classiques.

\begin{exemple}
La droite réelle.
\end{exemple}
Dans ce cas, $\rp1(M)=0$ et $\tp1(M)\not=0$. En effet, il existe des fonctions $L^p$ sur $\R$ dont les primitives ne sont pas $L^p$, mais toute fonction $L^p$ est limite $L^p$ de fonctions à support compact, d'intégrale nulle, qui admettent une primitive à support compact (donc $L^p$).

\begin{exemple}
Le plan hyperbolique.
\end{exemple}
Dans ce cas, $\tp1(M)=0$ et $\rp1(M)\not=0$. L'annulation de la torsion résulte de l'inégalité isopérimétrique linéaire : l'aire d'un domaine du plan hyperbolique est inférieure à la longueur de son bord. Lorsque $p=2$, la cohomologie réduite est représentée par des 1-formes harmoniques $L^2$. Leurs primitives possèdent des valeurs au bord, des fonctions sur le cercle qui appartiennent à l'espace de Sobolev $H^{1/2}$. Plus généralement, pour tout $p>1$, $\rp1(M)$ est isomorphe à l'espace de Besov $B^{p,p}_{1/p}$ du cercle, quotienté par les fonctions constantes, voir \cite{Triebel}.

\medskip

La droite réelle illustre le cas particulier des espaces homogènes \emph{fermés à l'infini}, i.e. pour lesquels l'inégalité isopérimétrique linéaire est en défaut, voir Définition \ref{ouverte}. Ces espaces ont été caractérisés par H. Hoke, \cite{Hoke}. Le plan hyperbolique illustre une autre famille particulière, celle des espaces homogènes à courbure sectionnelle négative. Ces espaces ont été caractérisés par E. Heintze, \cite{Heintze}. Le résultat principal de cet article est le fait que dans le cas général, i.e. avec pour seules exceptions ces deux familles, la cohomologie $L^p$ en degré 1 est nulle. En particulier, les espaces homogènes possédant des fonctions harmoniques non constantes dont le gradient est $L^p$ sont rares.

\subsection{Le résultat}

Remarquer que l'espace $\hp1(M)$ ne change pas si on change de métrique riemannienne sur $M$, pourvu qu'elle reste bornée par rapport à la métrique initiale (plus généralement, $\hp1(M)$ est un invariant de quasiisométrie, voir au paragraphe \ref{qi}). Par conséquent, si $M$ est homogène, $\hp1(M)$ ne change pas si on change de métrique riemannienne sur $M$, pourvu qu'elle reste invariante par un groupe cocompact d'isométries. Plus généralement, c'est un invariant de commensurabilité.

\begin{defi}
\label{commensurable}
On dit que deux groupes de Lie connexes $G$ et $G'$ sont commensurables s'il sont dans la même classe pour la relation d'équivalence engendrée par les relations suivantes.
\begin{itemize}
\item $G$ est isomorphe à un sous-groupe fermé, connexe et cocompact de $G'$.
\item $G'=G/K$ est le quotient de $G$ par un sous-groupe compact distingué.
\end{itemize}
On dit que deux espaces homogènes riemanniens $M$ et $M'$ sont \emph{commensurables} si les composantes neutres de leurs groupes d'isométries le sont.
\end{defi}

\begin{theo}
\label{degre1}
Soit $M$ un espace homog{\`e}ne riemannien non compact. Alors la cohomologie $L^p$ en degr{\'e} $1$ de $M$ est nulle, sauf si $M$ est ferm{\'e} {\`a} l'infini ou {\`a} courbure n{\'e}gative. Plus pr{\'e}cis{\'e}ment,
\begin{enumerate}
\item ou bien le groupe d'isom{\'e}tries de $M$ est une extension compacte d'un groupe de Lie r{\'e}soluble unimodulaire. Dans ce cas, pour tout $p>1$, $\tp1(M)\not=0$.
\item ou bien $M$ est commensurable {\`a} un espace homog{\`e}ne {\`a} courbure sectionnelle strictement n{\'e}gative. Dans ce cas, $\tp1(M)=0$ pour tout $p\geq 1$. De plus, il existe un nombre ${\bf p}(M)\geq 1$ tel que $\hp1(M)=0$ si $p\leq {\bf p}(M)$ et $\rp1(M)\not=0$ si $p>{\bf p}(M)$.
\item sinon, $\hp1(M)=0$ pour tout $p>1$.
\end{enumerate}
\end{theo}

Pour un espace homog\`ene \`a courbure strictement n\'egative, le nombre ${\bf p}(M)$, qui sera explicité dans la Proposition \ref{homneg}, peut s'interpr\'eter comme la dimension de Hausdorff du bord \`a l'infini, voir \cite{Pdim}. Il vaut au moins $2$, sauf pour le plan hyperbolique.

\begin{cor}
\label{car2}
Soit $M$ un espace homog{\`e}ne riemannien non compact. Si $M$ poss{\`e}de des fonctions harmoniques non constantes dont le gradient est $L^2$, il est commensurable (et donc quasiisom{\'e}trique) au plan hyperbolique.
\end{cor}

Dans la classe plus restreinte des espaces symétriques riemanniens, ce résultat est dû à A. Borel, \cite{Borel}.

\medskip

Un résultat récent de R. Tessera \cite{Tessera} affirme que, pour un groupe de Lie r{\'e}soluble unimodulaire, $\rp1(M)=0$. Cela permet de compléter le théorème \ref{degre1} comme suit. \emph{Pour un espace homogène riemannien, on a la trichotomie suivante.
\begin{enumerate}
\item ou bien, pour tout $p>1$, $\rp1(M)=0$ et $\tp1(M)\not=0$.
\item ou bien, pour tout $p>1$, $\tp1(M)=0$, et il existe un nombre ${\bf p}(M)\geq 1$ tel que $\rp1(M)\not=0\Leftrightarrow p>{\bf p}(M)$.
\item ou bien, pour tout $p>1$, $\hp1(M)=0$.
\end{enumerate}
}
\subsection{Schéma de la preuve}

\begin{enumerate}
\item Par commensurabilité, on se ramène au cas des groupes de Lie résolubles connexes $G$ munis de métriques riemanniennes invariantes à gauche.
\item Si $G$ est unimodulaire, alors, d'après un résultat de H. Hoke \cite{Hoke}, il est fermé à l'infini, donc, d'après un principe qui remonte à H. Kesten \cite{Kesten}, la torsion $L^p$ est non nulle pour tout $p\geq 1$.
\item Si $G$ n'est pas unimodulaire, on choisit un sous-groupe à un paramètre $t\mapsto g_t$ qui augmente le volume. En s'inspirant de R. Strichartz, \cite{Strichartz}, on montre que, pour toute fonction $u$ sur $G$ dont le gradient est $L^p$, la limite $u_{\infty}(x)=\lim_{t\to +\infty}u(xg_t)$ existe pour presque tout $x$, et que $u\in L^p$ si et seulement si $u_{\infty}$ est constante. Autrement dit, l'application $u\mapsto u_{\infty}$ induit un isomorphisme de l'espace $\hp1(G)$ sur le sous-espace, noté $\mathcal{B}^p$, des fonctions $g_t$-invariantes ayant $-1$ dérivées dans $L^{p}$, modulo les fonctions constantes. En particulier, $\tp1(G)=0$.
\item On montre que toute fonction de $\mathcal{B}^p$ est invariante à droite par les éléments de $G$ qui sont contractés (au sens large) par $Ad(g_{t})$. Or une fonction de norme $L^{p}_{-1}$ finie ne peut-être invariante à gauche par un sous-groupe distingué non compact de $G$. On conclut que si $\mathcal{B}^p \not=\{0\}$, $Ad(a_t )$ dilate dans toutes les directions, i.e. que $G$ possède des métriques riemanniennes invariantes à gauche à courbure sectionnelle strictement négative, \cite{Heintze}.
\end{enumerate}

\subsection{Questions}

\begin{enumerate}
\item Du théorème \ref{degre1}, il résulte que, parmi les espaces homogènes riemanniens, la classe des espaces homogènes commensurables à un espace à courbure négative est invariante par quasiisométries. Peut-on le voir de façon plus directe ? Par exemple, Y. de Cornulier conjecture que cette classe contient exactement les espaces homogènes riemanniens qui sont hyperboliques au sens de M. Gromov.

\item D'un résultat de N. Lohoué, \cite{Lohoue} (voir aussi \cite{DAR}), il résulte que, pour les espaces ouverts à l'infini, le fait que la cohomologie réduite soit non nulle, $\rp1(M)\not=0$, équivaut à l'existence de fonctions harmoniques non constantes dont le gradient est $L^p$. Le théorème \ref{degre1} répond donc à la question de l'existence de fonctions harmoniques non constantes dont le gradient est $L^p$ sur les espaces homogènes riemanniens ouverts à l'infini. Dans le cas fermé à l'infini, la question reste ouverte.
\end{enumerate}

\subsection{Remerciements}

Je  suis reconnaissant à Y. Benoist et à Y. de Cornulier pour leur aide au sujet de la structure des groupes de Lie, ainsi qu'à R. Tessera pour l'intérêt qu'il a porté à ce texte, extrait d'un manuscrit ancien non publié, \cite{P'}.

\section{Réduction aux groupes résolubles}

\subsection{Invariance sous quasiisométrie}
\label{qi}

Une application $f:X\to X'$ entre espaces m{\'e}triques est une {\em
quasiisom{\'e}trie} s'il existe des constantes $L$ et $D$ telles que tout point
de $X'$ est {\`a} distance au plus $D$ de l'image de $f$, et pour tous $x$,
$y\in X$, 
$$\d -D+\inv{L}d^{X}(x,y)\leq d^{X'}(f(x),f(y))\leq L\,d^{X}(x,y)+D.$$

\begin{lemme}
\label{cocompact}
Deux espaces homogènes riemanniens connexes commensurables sont qua\-si\-i\-so\-mé\-tri\-ques.
\end{lemme}

\preuve
Soit $G$ un groupe de Lie connexe, soit $H\subset G$ un sous-groupe compact. Montrons que les espaces homogènes $G$ et $M=G/H$ sont quasiisométriques.

Choisissons un suppl{\'e}mentaire $ad_H$-invariant $U$ de l'alg{\`e}bre de Lie ${\cal H}$ dans ${\cal G}$ et une structure euclidienne $ad_H$-invariante $g_{U}$ sur $U$. Cela d{\'e}termine une m{\'e}trique riemannienne $G$-invariante sur $M$. Soit $g_{{\cal H}}$ une structure euclidienne sur ${\cal H}$. La m{\'e}trique euclidienne $g_{{\cal H}}\oplus g_{U}$ sur ${\cal G}$ d{\'e}termine une m{\'e}trique riemannienne invariante {\`a} gauche sur $G$, invariante par l'action {\`a} droite de $H$. Pour cette m{\'e}trique, l'application
$f:G\mapsto G/H$ est une submersion riemannienne, i.e. isom{\'e}trique orthogonalement aux fibres. Par cons{\'e}quent, $f$ diminue les distances.

Inversement, {\'e}tant donn{\'e} un point $g\in G$, toute courbe dans $M=G/K$ d'origine $f(g)$ a un unique rel{\`e}vement en une courbe d'origine $g$ dans $G$ partout orthogonale aux fibres, et le rel{\`e}vement pr{\'e}l{\`e}ve la longueur. Par cons{\'e}quent, la distance entre deux fibres $f^{-1}(x)$ et $f^{-1}(x')$ est {\'e}gale {\`a} $d^{M}(x,x')$. Enfin, par invariance {\`a} gauche, toutes les fibres sont
isom{\'e}triques de m{\^e}me diam{\`e}tre $D$. Par cons{\'e}quent, pour tous $g$, $g'\in G$, $\d d^{M}(f(g),f(g'))\geq d^{G}(g,g')-D$. Ceci montre que $f$ est une quasiisom{\'e}trie, et $M$ est quasiisom{\'e}trique {\`a} $G$.

Lorsque $H$ est compact et distingué dans $G$, l'action de $G$ sur $G/H$ passe au quotient en l'action de $G'=G/H$ sur lui-même par translation à gauche, donc les métriques $G$-invariantes sur $G/H$ coïncident avec les métriques invariantes à gauche sur $G'$. 
On conclut que $G$ et $G'$ sont quasiisométriques.

Soit maintenant $G$ un groupe de Lie connexe et $H$ un sous-groupe fermé connexe cocompact de $G$. Montrons que $G$ et $H$ (munis de métriques riemanniennes invariantes à gauche) sont quasiisométriques.

Montrons d'abord par l'absurde qu'il existe une constante $D$ telle que  $d^{G}(g,H)\leq D$ pour tout $g\in G$. Sinon, il existe une suite $g_j$ telle que $d^{G}(g_{j},H)$ tend vers $+\infty$. L'espace des classes {\`a} gauche $H\setminus G$ est lui aussi compact. Il existe donc une suite $h_{j}\in H$ telle que la suite $h_{j}^{-1}g_j$ converge dans $G$. Les distances
$d^{G}(h_{j},g_{j})=d^{G}(h_{j}^{-1}g_{j},1)$ sont born{\'e}es, et 
$d^{G}(g_{j},H)\leq d^{G}(h_{j},g_{j})$ l'est aussi, contradiction.

Soient $h$, $h'\in H$.  Alors $\delta=d^{G}(h,h')\leq d^{H}(h,h')$. Il reste {\`a} montrer l'in{\'e}galit{\'e} inverse. Soit $\gamma:[0,\delta]\to G$ une g{\'e}od{\'e}sique minimisante de $h$ {\`a} $h'$. Pour $i=0,\ldots,N=[\delta/D]$, posons $g_{i}=\gamma(iD)$, de sorte que $d^{G}(g_{i},g_{i+1})\leq D$ et $ND\leq\delta\leq
(N+1)D$. Choisissons $h_{i}\in H$ tel que $d^{G}(h_{i},g_{i})\leq D$. Comme $H$ est ferm{\'e} dans $G$, l'intersection de $H$ et de la boule de rayon $3D$ de $G$ est compacte, donc contenue dans une boule de rayon $D'$ pour la m{\'e}trique de $H$. Pour tout $i$, $d^{G}(h_{i}^{-1}h_{i+1},1)\leq 3D$, donc $d^{H}(h_{i}^{-1}h_{i+1},1)\leq D'$. Par cons{\'e}quent,
$$
d^{H}(h,h')\leq \sum d^{H}(h_{i}^{-1}h_{i+1},1)\leq (N+1)D'
\leq \frac{D'}{D} d^{G}(h,h')+D.\qed
$$

\begin{prop}
\label{exacte}
Soit $M$ un espace homog{\`e}ne riemannien connexe et non compact. Alors toute forme ferm{\'e}e $L^p$ sur $M$ est exacte.
\end{prop}

\preuve
Soit $M$ un espace homog{\`e}ne riemannien connexe et non compact. La composante neutre $G$ du groupe d'isom{\'e}tries de $M$ est transitive sur $M$. Soit $\omega$ une forme $L^p$ ferm{\'e}e non exacte sur $M$. Alors il existe une boule $B$ de $M$ telle que la restriction de $\omega$ {\`a} $B$ ne soit pas exacte.

Lorsque $g$ varie dans $G$, les formes $g^{*}\omega$ sont deux {\`a} deux cohomologues. Par cons{\'e}quent la classe de cohomologie $c\in H^{k}(B,\R)$ de la restriction de $g^{*}\omega$ {\`a} $B$ est ind{\'e}pendante de $g\in G$. Comme $B$ est compacte, la cohomologie $L^p$ de $B$ coïncide avec sa cohomologie de de Rham ordinaire. En particulier, la différentielle extérieure $d$ agissant sur les formes différentielles $L^p$ (et à différentielle $L^p$) sur $B$ a une image fermée, et il existe une constante $\epsilon>0$ telle que toute forme ferm{\'e}e repr{\'e}sentant la classe $c$ soit de norme $L^p$ sup{\'e}rieure {\`a} $\epsilon$. Par cons{\'e}quent $\d \n{\omega}_{L^{p}(gB)}=\n{g^{*}\omega}_{L^{p}(B)}\geq\epsilon$. C'est incompatible avec l'hypoth{\`e}se $\omega\in L^{p}(M)$.\qed

\begin{cor}
\label{h1qi}
$\hp1$, $\tp1$ et $\rp1$ sont des invariants de commensurabilité pour les espaces homog{\`e}nes riemanniens connexes non compacts.
\end{cor}

\preuve
Il est classique (\cite{Gromovasympt,BourdonPajot}) que, dans la catégorie des variétés riemanniennes à géométrie bornée, la cohomologie $L^p$ \emph{exacte} en degré 1 est un invariant de quasiisométrie (l'énoncé précis se trouve dans \cite{P"}). C'est donc un invariant de commensurabilité.\qed

\subsection{Réduction aux groupes résolubles}

\begin{lemme}
\label{qires}
Notons $T_n$ le groupe des matrices triangulaires réelles de taille $n$. Tout espace homog{\`e}ne riemannien connexe est commensurable {\`a} un sous-groupe fermé connexe de $T_n$, pour un $n$ assez grand.
\end{lemme}

\preuve
Soit $M=\breve{G}/\breve{H}$ un espace homog{\`e}ne riemannien connexe. Alors les orbites dans $M$ de l'action de la composante neutre $G$ de $\breve{G}$ sont ouvertes, donc elles sont fermées. Comme $M$ est connexe, il n'y en a qu'une, et on peut écrire $M=G/H$ où $H=\breve{H}\cap G$ est compact. Par définition, $M$ est commensurable {\`a} $G$.

Soit $G=RS$ la d{\'e}composition de Levi du groupe $G$, o{\`u} $R$ est le radical de $G$ et $S$ est semi-simple (\cite{Bourbaki} chap. 3 page 244). Supposons d'abord que le centre de $S$ est fini. Soit $S=NAK$ la d{\'e}composition d'Iwasawa de $S$. Comme le centre de $S$ est fini, le groupe $K$ est compact. Alors $G'=RNA$ est r{\'e}soluble et connexe, $G/G'=G/RNA=S/NA$ est compact, donc $G$ et $G'$ sont commensurables. 

Dans le cas général, une construction supplémentaire permet de se débarrasser du centre de $S$. L'image de $S$ dans le groupe des automorphismes de l'algèbre de Lie de $R$ est un groupe de Lie semi-simple connexe linéaire, donc son centre est fini. Par conséquent, il existe un sous-groupe ab{\'e}lien libre $\Gamma$ d'indice fini dans le centre de $S$ qui agit trivialement sur $R$. $\Gamma$ est donc contenu dans le centre de $G$. Posons $C=\Gamma\otimes\R$ et $G''=(G\times C)/\Gamma$ o{\`u} $\Gamma$ agit diagonalement par translation. Comme $G\times C$ est connexe, $G''$ l'est aussi. Comme $G''/G$, homéomorphe à $C/\Gamma$, est compact, $G$ est cocompact dans $G''$, donc $G$ est commensurable {\`a} $G''$. 

La projection sur le premier facteur passe au quotient en $\pi:G''\to \overline{G}=:G/\Gamma$. Le noyau de $\pi$ est isomorphe {\`a} $C$ donc connexe. Soit $\overline{S}=S/\Gamma$ et $\overline{R}=R/\Gamma\cap R$. Ces sous-groupes de $\overline{G}$ engendrent $\overline{G}$. Soit $\overline{S}=NAK$ la d{\'e}composition d'Iwasawa de $\overline{S}$. A nouveau, comme le centre de $\overline{S}$ est fini, le groupe $K$ est compact. Alors $G'=\pi^{-1}(RNA)$ est r{\'e}soluble et connexe, $G''/G'=\overline{G}/\overline{R}NA=\overline{S}/NA$ est compact, donc $G''$ (et donc $G$) est commensurable {\`a} $G'$. 


D'après le Lemme 2.4 de \cite{dC}, $G'$ est isomorphe à un sous-groupe fermé cocompact d'un groupe $H$ qui possède un sous-groupe fermé connexe cocompact $T$ qui se plonge proprement dans l'un des groupes $T_n$. Autrement dit, $G'$ est commensurable à l'image de $T$ dans $T_n$, qui est fermée.\qed

\section{Torsion et inégalité isopérimétrique}

\begin{prop}
\label{torsion1}
Soit $M$ une vari{\'e}t{\'e} riemannienne de volume infini. Soit
$p\geq 1$. La torsion $\tp1(M)$ est nulle si et seulement si l'in{\'e}galit{\'e} de Sobolev suivante est vraie. Il existe une constante telle que, pour toute fonction $u\in L^p (M)$ telle que $du\in L^p (M)$, 
$$
\n{u}\lp \leq\con\n{du}\lp.
$$
\end{prop}

\preuve
Notons $\Omega^{0,p} (M)$ l'ensemble des fonctions $L^p$ à gradient $L^p$, muni de la norme $\n{u}_{\Omega^{0,p}}=\n{u}\lp+\n{du}\lp$. Notons $\Omega^{1,p}(M)$ l'espace des 1-formes $L^p$ sur $M$. Ce sont des espaces de Banach, et $d:\Omega^{0,p} \to  \Omega^{1,p}$ est borné.

Supposons l'in{\'e}galit{\'e} de Sobolev vraie. Soit $\alpha$ une 1-forme $L^p$ dont la classe de cohomologie appartient à $\tp1(M)$. Soit $u_j$ une suite de fonctions $L^p$ sur $M$ telle que $du_j \in L^p $ et $du_j$ converge vers $\alpha$ dans $L^p$. D'apr{\`e}s l'in{\'e}galit{\'e} de Sobolev, la suite $u_j$ est born{\'e}e dans $L^p$. Extrayons une sous-suite qui converge faiblement vers $u_\infty$. Alors $u_\infty\in L^p$ et $du_{\infty}=\alpha$, donc la classe de cohomologie $L^p$ de $\alpha$ est nulle. On conclut que $\tp1(M)=0$.

Inversement, supposons que $\tp1(M)=0$. Alors $d\Omega^{0,p}(M)$ est fermé dans $\Omega^{1,p}(M)$, donc c'est un espace de Banach. Comme $M$ est de volume infini, les fonctions constantes ne sont pas $L^p$, donc $d$ est injectif. L'op{\'e}rateur $d$ induit $\overline{d}:\Omega^{0,p}(M)\to d\Omega^{0,p}(M)$. C'est une bijection continue entre espaces de Banach, donc un isomorphisme. Par conséquent, il existe une constante telle que pour tout $u\in \Omega^{0,p}(M)$, $\n{u}_{\Omega^{0,p}}\leq\con\n{du}\lp$.\qed

\begin{defi}
\label{ouverte}
Une vari{\'e}t{\'e} riemannienne $M$ est dite {\em ouverte {\`a} l'infini} si elle satisfait une in{\'e}galit{\'e} isop{\'e}rim{\'e}trique lin{\'e}aire, i.e. si pour toute sous-vari{\'e}t{\'e} compacte $\Delta$ de codimension $0$, {\`a} bord lisse,
$$
{\rm Vol}\,\Delta\leq\con{\rm Vol}\,\partial\Delta.
$$
Si $M$ n'est pas ouverte {\`a} l'infini, elle est {\em ferm{\'e}e {\`a} l'infini}.
\end{defi}

\begin{prop}
Soit $M$ une vari{\'e}t{\'e} riemannienne {\`a} g{\'e}om{\'e}trie born{\'e}e. Les propri{\'e}t{\'e}s suivantes sont {\'e}quivalentes.
\begin{enumerate}
\item $M$ est ouverte {\`a} l'infini ;
\item pour tout $p\geq 1$, $\tp1(M)=0$ ;
\item il existe $p\geq 1$ tel que $\tp1(M)=0$.
\end{enumerate}
\end{prop}

\preuve
Lorsque $p=2$, ces {\'e}quivalences sont essentiellement dues {\`a} R. Brooks, \cite{Brooks}. Elles apparaissent pour la première fois dans le contexte des marches aléatoires sur les groupes, \cite{Kesten}.

Il est classique (voir \cite{Mazya}) que l'in{\'e}galit{\'e} isop{\'e}rim{\'e}trique ${\rm Vol}\,\Delta\leq\con{\rm Vol}\,\partial\Delta$ est {\'e}quivalente {\`a} l'in{\'e}galit{\'e} de Sobolev $L^1$ (avec la m{\^e}me constante) : pour toute fonction $u$ lisse {\`a} support compact sur $M$, $\d\n{\phi}_{L^{1}}\leq\con\n{du}_{L^{1}}$. 

L'in{\'e}galit{\'e} de Sobolev $L^1$, appliqu{\'e}e {\`a} la fonction $u^p$, combin{\'e}e {\`a} l'in{\'e}galit{\'e} de H{\"o}lder, entra{\^\i}ne l'in{\'e}galit{\'e} de Sobolev $L^p$ pour $p\geq 1$. Ceci prouve que $(1)\Rightarrow(2)$.

Il reste {\`a} montrer que, sous l'hypoth{\`e}se de g{\'e}om{\'e}trie born{\'e}e, l'i\-n{\'e}\-ga\-li\-t{\'e} de Sobolev $L^p$ pour un $p>1$ entra{\^\i}ne l'in{\'e}galit{\'e} de Sobolev $L^1$. Par invariance par quasiisom{\'e}trie, on peut remplacer $M$ par un graphe $X$ de valence born{\'e}e par $V$ et dont les ar{\^e}tes ont toutes m{\^e}me longueur. Les fonctions sur $M$ deviennent des fonctions sur l'ensemble $X^0$ des sommets de $X$ et la norme $L^p$ de la diff{\'e}rentielle devient $\d(\sum|u(a)-u(b)|^{p})^{1/p}$ o{\`u}
la somme est {\'e}tendue aux couples de sommets reli{\'e}s par une ar{\^e}te.

Supposons d'abord que $u$ ne prend que les valeurs $0$ et $1$. Alors $\d
\n{u}_{L^{p}}=\n{u}_{L^{1}}^{1/p}$ et $\d \n{du}_{L^{1}}=\n{du}_{L^{1}}^{1/p}$ donc l'in{\'e}galit{\'e} de Sobolev $L^p$ entra{\^\i}ne l'in{\'e}galit{\'e} de Sobolev $L^1$ dans ce cas. 

Enfin, on {\'e}crit une fonction quelconque comme int{\'e}grale de fonctions {\`a} valeurs dans $\{0,1\}$. Soit $u$ une fonction sur
$X^0$. Quitte {\`a} remplacer $u$ par sa valeur absolue (ce qui ne change pas la norme $L^1$ et diminue la norme $L^1$ de la diff{\'e}rentielle), on peut supposer que $u$ est positive ou nulle. On {\'e}crit alors $\d u=\int_{0}^{+\infty} u_{t}\,dt$ o{\`u} $\d
u_{t}(x)=1$ si $u(x)\geq t$ et $\d u_{t}(x)=0$ sinon. Alors
$$
\n{u}_{L^{1}}=\int_{0}^{+\infty} \n{u_{t}}_{L^{1}}\,dt
$$
et
$$
\n{du}_{L^{1}}=\int_{0}^{+\infty} \n{du_{t}}_{L^{1}}\,dt,
$$
donc l'in{\'e}galit{\'e} de Sobolev $L^p$ entra{\^\i}ne l'in{\'e}galit{\'e} de Sobolev $L^1$ en g{\'e}n{\'e}ral.\qed

\begin{prop}
\label{hoke}
{\em (H. Hoke, \cite{Hoke})}. Soit $G$ un groupe de Lie muni d'une m{\'e}trique riemannienne invariante {\`a} gauche. Alors $G$ est ferm{\'e} {\`a} l'infini si et seulement si $G$ est moyennable (i.e. extension d'un groupe résoluble par un groupe compact) et unimodulaire.
\end{prop}

\begin{cor}
\label{unimod}
Soit $G$ un groupe de Lie. Si $G$ est moyennable unimodulaire, alors pour tout $p>1$, $\tp1(G)\not=0$. Inversement, s'il existe $p>1$ tel que $\tp1(G)\not=0$, alors $G$ est moyennable unimodulaire.
\end{cor}

\section{Valeur au bord}

Soit $G$ un groupe de Lie non unimodulaire muni d'une métrique invariante à gauche, soit $\G$ son algèbre de Lie. La forme lin{\'e}aire $\xi\mapsto\tr ad_{\xi}$ sur $\G$ n'est pas nulle. Choisissons $\xi_0 \in\G$ tel que $\tau=\tr ad_{\xi_0}>0$. Notons $\phi_t =R_{\exp(t\xi_0 )}$ la translation à droite par $\exp(t\xi_0 )$. C'est le flot de $\xi_0$ vu comme champ de vecteurs invariant à gauche sur $G$, donc il multiplie les volumes par $e^{t\tau}$. Alors, pour toute fonction $u$ sur $G$ dont le gradient est $L^p$ et tout $t\in\R$,
\begin{eqnarray*}
\left\|\frac{d}{dt}u\circ\phi_t\right\|\lp &=&\n{(\xi_0 u)\circ\phi_t}\lp\\
&=&e^{-t\tau/p}\n{\xi_0 u}\lp\\
&\leq&\con\, e^{-t\tau/p}\n{du}\lp,
\end{eqnarray*}
où la constante est la norme (constante) de $\xi_0$. Comme $\tau/p>0$, l'application $t\mapsto u\circ\phi_t -u$ converge dans $L^p (G)$ quand $t$ tend vers $+\infty$. En particulier, $u\circ\phi_t$ converge presque partout.

\begin{defi}
\label{uinfinity}
Soit $u$ une fonction sur $G$ dont le gradient est $L^p$. On appelle \emph{valeur au bord} de $u$, et on note $u_{\infty}$ la limite des translatées $u\circ R_{\exp(t\xi_0 )}=u\circ\phi_t$, quand $t$ tend vers $+\infty$.
\end{defi}

Le terme valeur au bord s'inspire du cas de la courbure strictement négative. Lorsque $G$ possède une métrique riemannienne invariante à gauche à courbure sectionnelle négative, les orbites de $\xi_0$ convergent en $-\infty$ vers un même point du bord à l'infini, et en $+\infty$ vers tous les autres points du bord. 

On va montrer qu'une fonction sur $G$ dont le gradient est $L^p$ satisfait l'inégalité de Sobolev dès que sa valeur au bord est presque partout nulle. Cette idée est due à R. Strichartz, \cite{Strichartz}, voir aussi \cite{GKS}.

\subsection{Les espaces $(L^p +dL^p )(G)$ et $d^{-1}dL^{p}(G)$}

Dans cette section, on précise la régularité de la valeur au bord. Soit $u$ une fonction sur $G$ dont le gradient est $L^p$. Par construction, $v=u_{\infty}-u\in L^p$, donc $du_{\infty}$ est la somme d'une 1-forme $L^p$ et de la différentielle d'une fonction $L^p$. Cela suggère la définition suivante.

\begin{defi}
\label{lp+dlp}
Soit $\alpha$ une 1-forme (à coefficients distributions) sur $G$, on pose
\begin{eqnarray*}
\n{\alpha}_{(L^p +dL^{p})(G)}=\inf\{\n{\beta}\lp +\n{\gamma}\lp \,|\,\beta \textrm{ 1-forme, }\gamma \textrm{ fonction, }\alpha=\beta+d\gamma\}.
\end{eqnarray*}
\end{defi}

On va exprimer cette norme sous une forme plus commode, à l'aide de la définition suivante.

\begin{defi}
\label{lp-1}
Soit $B$ une boule de rayon 1 de $G$. Soit $f$ une fonction sur $B$. On pose
\begin{eqnarray*}
\n{f}_{d^{-1}dL^{p}(B)}=\inf_{c\in\R}\n{f-c}_{L^p (B)}.
\end{eqnarray*}

Pour une fonction $f$ localement $L^p$ sur $G$, on pose
\begin{eqnarray*}
\n{f}_{d^{-1}dL^{p}(G)}^{p}=\int_{G}\n{f_{|gB}}_{d^{-1}dL^{p}(gB)}^{p}\,dg.
\end{eqnarray*}
\end{defi}

Autrement dit, une fonction est de classe $d^{-1}dL^{p}$ si sa dérivée est localement la dérivée d'une fonction $L^p$, et si la fonction ``meilleure norme $L^p$ d'une primitive locale'' est elle-même $L^p$ sur $G$ tout entier.

\begin{lemme}
\label{d-1d}
Soit $w$ une fonction sur $G$. Alors $w$ appartient à $d^{-1}dL^{p}(G)$ si et seulement si $dw$ appartient à $(L^p +dL^{p})(G)$.
\end{lemme}

\preuve
Si $f$ est une fonction $L^p$ sur $G$, alors pour toute boule de rayon 1, $\n{f_{|B}}_{d^{-1}dL^{p}(B)}\leq \n{f_{|B}}_{L^{p}(B)}$, d'où, en intégrant sur toutes les boules, $\n{f}_{d^{-1}dL^{p}(G)}\leq\n{f}\lp$.

Soit $f$ une fonction sur $G$ dont la différentielle est $L^p$. D'après l'inégalité de Poincaré, il existe pour chaque boule $B$ une constante $c_{f,B}$ telle que
\begin{eqnarray*}
\n{f-c_{f,B}}_{L^p (B)}\leq\con\n{df}_{L^p (B)}.
\end{eqnarray*}
En intégrant sur toutes les boules de rayon 1, il vient $\n{f}_{d^{-1}dL^{p}(G)}\leq\con\n{df}\lp$.

Soit $w$ une fonction telle que $dw=\beta+d\gamma\in (L^p +dL^{p})(G)$. Alors $w-\gamma$ a une différentielle $L^p$, donc elle appartient à $d^{-1}dL^{p}(G)$, d'où $w=w-\gamma+\gamma \in d^{-1}dL^{p}(G)$. 

Réciproquement, si $w\in d^{-1}dL^{p}(G)$, construisons une fonction lisse, qui diffère de $w$ d'une fonction $L^p$, et dont la différentielle est $L^p$. Pour cela, il suffit de régulariser $w$. Soit $\psi$ une fonction lisse, positive ou nulle, à support compact dans la boule unité $B$ de $G$, d'intégrale 1. Soit $u$ le produit de convolution $u=w\star \psi$. Alors $u$ est lisse. Comme, pour toute constante $c$,  $du=(w-c)\star d\psi$ est une combinaison linéaire de valeurs de $w-c$, pour tout $g\in G$,
\begin{eqnarray*}
\n{du}_{L^p (gB)}\leq\con\inf_{c\in\R} \n{w-c}_{L^p (gB)}\leq\con \n{dw}_{d^{-1}dL^{p}(gB)},
\end{eqnarray*}
d'où, en intégrant par rapport à $g$, $du\in L^p (G)$. De même, pour toute constante $c$, $u-w=(w-c)\star(\psi-\delta)$, où $\delta$ est la masse de Dirac à l'origine, donc pour tout $g\in G$, 
\begin{eqnarray*}
\n{u-w}_{L^p (gB)}\leq\con\inf_c \n{w-c}_{L^p (gB)}\leq\con \n{dw}_{d^{-1}dL^{p}(gB)},
\end{eqnarray*}
d'où $u-w\in L^p (G)$. Cela prouve que $dw\in (L^p +dL^p )(G)$.\qed

\subsection{Image de l'application valeur au bord}

\begin{defi}
\label{besov}
Soit $G$ un groupe de Lie non unimodulaire. Soit $\xi_0 \in\G$ un champ de vecteurs invariant à gauche de divergence strictement positive. On note
\begin{eqnarray*}
\mathcal{B}_{\xi_0}^{p} =\{w\in d^{-1}dL^{p}(G)\,|\,w\circ R_{\exp(t\xi_0 )}=w \textrm{ presque partout, pour tout }t\},
\end{eqnarray*}
\end{defi}

\begin{prop}
\label{strichartz}
Soit $G$ un groupe de Lie non unimodulaire. On choisit un vecteur $\xi_0 \in\G$ de divergence positive, pour définir des valeurs au bord. Alors l'application $u\mapsto u_{\infty}$ passe au quotient en une injection
$\hp1(G)\to d^{-1}dL^{p}(G)$ mod $\R$, dont l'image est exactement $\mathcal{B}_{\xi_0}^{p}$ mod $\R$, i.e. modulo les constantes additives. 

De façon équivalente, l'application $du\mapsto du_{\infty}$ passe au quotient en une injection $\hp1(G)\to (L^p +dL^{p})(G)$, dont l'image est le sous-espace
\begin{eqnarray*}
d\mathcal{B}_{\xi_0}^{p} =\{\alpha\in (L^p +dL^{p})(G)\,|\,d\alpha=0 \textrm{ et }\alpha(\xi_0 )=0 \textrm{ au sens des distributions}\}.
\end{eqnarray*}
\end{prop}

\preuve
Les valeurs au bord $u_\infty$ sont des fonctions $\phi_t$-invariantes. Par construction, leur diffé\-ren\-ti\-el\-le appartient à $L^p +dL^{p}$. D'après le lemme \ref{d-1d}, $u_\infty \in d^{-1}dL^{p}$, donc $u\in \mathcal{B}_{\xi_0}^{p}$. 

Si $u\in L^p$, $\n{u\circ\phi_t}\lp =e^{-t\tau/p}\n{u}\lp$ tend vers 0 quand $t$ tend vers $+\infty$, donc $u_{\infty}=0$. Si $u$ est constante, $u_{\infty}$ l'est aussi. Par conséquent, l'opérateur valeur au bord, à constante additive près, passe au quotient. Si $u_\infty$ est constante, comme $u-u_\infty \in L^p$, la classe de cohomologie de $du$ est nulle. Par conséquent, l'opérateur valeur au bord est injectif.

Soit $\alpha$ une 1-forme fermée de $(L^p +dL^p )(G)$, telle que $\alpha(\xi_0 )=0$ au sens des distributions. Par hypothèse, il existe une 1-forme $L^p$ $\beta$ et une fonction $L^p$ $\gamma$ telles que $\alpha=\beta+d\gamma$. Nécessairement, $d\beta=0$, donc, avec le lemme \ref{exacte}, $\beta=du$ où la fonction $u$ a un gradient $L^p$. Soit $w=u+\gamma$. Alors $dw=\alpha$, sa dérivée au sens des distributions le long des orbites de $\xi_0$ s'annule, donc $w$ est $\phi_t$-invariante. Comme $\gamma$ est $L^p$, les fonctions $\gamma\circ\phi_t$ tendent vers 0 dans $L^p$ quand $t$ tend vers $+\infty$, donc $w$ est la valeur au bord de $u$. On conclut que l'image de l'opérateur valeur au bord est exactement $\mathcal{B}_{\xi_0}^{p}$ modulo les constantes.\qed

\subsection{Restrictions sur les valeurs au bord}

Une métrique invariante à gauche sur un groupe de Lie $G$ croît exponentiellement avec des taux différents suivant les directions. Cela entraîne des restrictions sur les fonctions de $d^{-1}dL^{p}(G)$ invariantes à droite par un sous-groupe à un paramètre.

\begin{lemme}
\label{restr}
Soit $G$ un groupe de Lie non unimodulaire, soit $\xi_0 \in\G$ de divergence $\tr ad_{\xi_0}$ strictement positive. Pour toute fonction $w\in\mathcal{B}_{\xi_0}^{p}$, la différentielle au sens des distributions $dw$ s'annule sur tous les sous-espaces caractéristiques de $ad_{\xi_0}$ relatifs à des valeurs propres de parties réelles strictement inférieures à $\d \frac{\tr ad_{\xi_0}}{p}$.
\end{lemme}

\preuve
On écrit $w=u+v$ où $v\in L^p (G)$ et $du\in L^p$. Comme $w$ est $\phi_t$ invariante, $w=u\circ\phi_t +v\circ\phi_t$. Or $v\circ \phi_t$ tend vers 0 quand $t$ tend vers $+\infty$. Par conséquent,
\begin{eqnarray*}
dw=\lim_{t\to+\infty} d(u\circ\phi_t )=\lim_{t\to+\infty} \phi_{t}^{*}du
\end{eqnarray*}
au sens des distributions.

On voit $du$ comme une fonction $G\mapsto \G^*$ et les éléments de $\G$ comme des champs de vecteurs invariants à gauche. Soit $\eta\in\G\otimes\C$. Alors, en tout point de $G$,
\begin{eqnarray*}
\langle\phi_{t}^{*}du,\eta\rangle&=&\langle du\circ\phi_t ,ad_{\exp(t\xi_0 )}(\eta)\rangle
\end{eqnarray*}
Soit $\lambda\in\C$ une valeur propre de $ad_{\xi_0}$. Supposons que $\eta\in\ker(ad_{\xi_0}-\lambda\,I)^{\ell}\subset \mathcal{G}$ appartient à l'espace caractéristique relatif à $\lambda$. Alors
$$
\exp(t\,ad_{\xi_0})(\eta)=e^{t\lambda}P(t),
$$
où $P$ est un polynôme en $t$ à valeurs dans $\G$. Comme les fonctions $du\circ\phi_t$ tendent vers 0 exponentiellement vite dans $L^p$, $\n{du\circ\phi_t}\lp\leq e^{-t\,\tr ad_{\xi_0}/p}$, il vient
\begin{eqnarray*}
\n{\langle\phi_{t}^{*}du,\eta\rangle}\lp&=&\n{e^{t\lambda}\langle du\circ\phi_t ,P(t)\rangle}\lp\\
&\leq&e^{t(\Re e(\lambda)-\tr ad_{\xi_0}/p)}|P(t)|.
\end{eqnarray*}
Si $\Re e(\lambda)<\tr ad_{\xi_0}/p$, $\langle\phi_{t}^{*}du,\eta\rangle$
tend vers 0 quand $t$ tend vers $+\infty$, et on conclut que $\langle dw,\eta\rangle=0$. \qed

\bigskip

On aura besoin aussi du cas limite du lemme \ref{restr}.

\begin{lemme}
\label{restr'}
Soit $G$ un groupe de Lie connexe. Soit $\xi\in\G$ un champ de vecteurs invariant à gauche. On suppose que le sous-groupe à un paramètre $A=\{\exp(t\xi)\,|\,t\in\R\}$ est fermé dans $G$. Soit $w\in d^{-1}dL^{p}(G)$ une fonction $A$-invariante à droite. Alors la différentielle au sens des distributions $dw$ s'annule sur les sous-espaces propres de $ad_{\xi_0}$ relatifs aux valeurs propres de partie réelle égale à $\d \frac{\tr ad_{\xi}}{p}$. En particulier, si $\tr ad_{\xi}=0$ et si $[\xi,\eta]=0$, alors $\langle dw,\eta\rangle=0$.
\end{lemme}

\preuve
Comme $A$ est fermé dans $G$, $\phi_{t}=R_{\exp(t\xi)}$ définit une action propre de $\R$ sur $G$.

Soit $w\in d^{-1}dL^{p}(G)$ une fonction $A$-invariante à droite. De nouveau, on l'écrit $w=u+v$ avec $v\in L^{p}$ et $du\in L^{p}$. Dans le cas où $\tr ad_{\xi}=0$, il n'est plus vrai que les $v\circ\phi_{t}$ tendent vers 0 dans $L^{p}$, mais on peut quand même affirmer qu'il existe une sous-suite $t_{j}$ tendant vers $+\infty$ telle que $v\circ\phi_{t_{j}}$ tend vers 0 au sens des distributions. En effet, dans ce cas, $\phi_{t}$ préserve la mesure et agit proprement sur $G$. Pour tout compact $K$ de $G$, il existe une sous-suite $t_{j}$ tendant vers $+\infty$ telle que les $\phi_{t_{j}}(K)$ soient deux à deux disjoints. Comme la suite $(\n{v}_{L^{p}(\phi_{t_{j}}(K))})_{j}$ est $\ell^{p}$, elle tend vers 0. Cela prouve que les $v\circ\phi_{t_{j}}$ tendent vers 0 dans $L^{p}(K)$.

Une variante de ce raisonnement s'applique à la fonction $v'=|\langle du,\eta\rangle|$, lorsque $\eta$ est un vecteur propre de $ad_{\xi}$ relatif à une valeur propre $\lambda$ telle que $\Re e(\lambda)=\tr ad_{\xi}/p$. Comme $v'\in L^{p}$, on peut supposer que les intégrales $\int_{\phi_{t_{j}}(K)}v'^{p}$ tendent vers 0. Or
\begin{eqnarray*}
\int_{\phi_{t_{j}}(K)}v'^{p}=\int_{K}|\langle \phi_{t_{j}}^{*}du,\eta\rangle|^{p}.
\end{eqnarray*}
Par conséquent, $\langle \phi_{t_{j}}^{*}du,\eta\rangle$ tend vers 0 au sens des distributions. Cela suffit pour conclure que $\langle dw,\eta\rangle=0$.\qed

\bigskip

Une fonction $L^{p}$ sur un groupe de Lie $G$ ne peut être invariante par un sous-groupe unimodulaire non compact. Il en est de même pour une fonction de $d^{-1}dL^{p}(G)$.

\begin{lemme}
\label{invariante}
Soit $G$ un groupe de Lie. Soit $w\in d^{-1}dL^{p}(G)$. Si $w$ est invariante par un sous-groupe fermé distingué unimodulaire non compact de $G$, alors $w$ est constante. 
\end{lemme}

\preuve
Comme $H$ est distingué, toute fonction $H$-invariante à droite est aussi $H$-ionvariante à gauche. Comme $H$ est fermé, distingué et unimodulaire, la mesure de Haar de $G$ se désintègre le long des orbites de $H$, par rapport à une mesure de Haar de $G/H$, \cite{Weil}. Si $w$ est $H$-invariante, la fonction $g\mapsto I(g)=\n{w}_{d^{-1}dL^p (gB)}^p$ est aussi $H$-invariante, elle passe au quotient en une fonction $J$ sur $G/H$, et son intégrale sur $G$ vaut
\begin{eqnarray*}
\n{w}_{d^{-1}dL^p (G)}^p =vol(H)\int_{G/H}J(\bar{g}) \,d\bar{g}.
\end{eqnarray*}
Si $H$ est non compact, son volume est infini, donc $J=0$, et $w$ est constante.\qed

\begin{cor}
\label{invariante'}
Soit $G$ un groupe de Lie, $N\subset G$ un sous-groupe fermé nilpotent connexe, simplement connexe, distingué, non trivial. Soit $w\in d^{-1}dL^{p}(G)$, soit $\xi\in\N$ un champ de vecteurs invariant à gauche non nul. Si, au sens des distributions, $\langle dw,\xi\rangle=0$, alors $w$ est constante. 
\end{cor}

\preuve
Notons $A$ le sous-groupe à un paramètre engendré par $\xi$. Comme $N$ est nilpotent connexe et simplement connexe, l'exponentielle de $N$ est un difféomorphisme, donc $A$ est fermé dans $N$ et donc dans $G$. D'autre part, $ad_{\xi}$ est nilpotent donc $\tr ad_{\xi}=0$.

Soit $\eta$ un vecteur non nul du centre de $\N$. Alors $[\xi,\eta]=0$. Le lemme \ref{restr'} donne que $\langle dw,\eta\rangle=0$. Par conséquent, le centre de $N$, qui est non compact, fermé et distingué dans $G$, préserve $w$. Le lemme \ref{invariante} s'applique, on conclut que $w$ est constante.\qed

\subsection{Cohomologie $L^p$ en degr{\'e} $1$ des groupes de Lie r{\'e}solubles}

\begin{prop}
\label{nonuni}
Soit $G$ un sous-groupe connexe et fermé du groupe des matrices triangulaires réelles. On suppose que $G$ n'est pas unimodulaire et qu'il n'admet pas de m{\'e}trique invariante {\`a} gauche {\`a} courbure sectionnelle strictement n{\'e}gative. Alors pour tout $p>1$, $\hp1(G)=0$.
\end{prop}

\preuve
D'apr{\`e}s \cite{Bourbaki}, chap. 7, pages 19-20, il existe deux sous-groupes nilpotents connexes $A$ et $N\subset G$ tels que $G=AN$, $N$ est fermé, distingu{\'e} et $[A,A]\subset N$. Notons $\A$ et $\N$ leurs alg{\`e}bres de Lie. L'alg{\`e}bre de Lie $\G$ se d{\'e}compose en sous-espaces caract{\'e}ristiques pour l'action de l'alg{\`e}bre de Lie nilpotente $\A$,
$$
\G\otimes\C=\bigoplus_{\lambda\in\Sigma} \G_{\lambda}
$$
o{\`u} les $\lambda$ sont des formes lin{\'e}aires sur $\A$ telles que pour $\xi\in\A$, $ad_{\xi}-\lambda(\xi)I$ restreinte {\`a} $\G_{\lambda}$ soit nilpotente (voir \cite{Bourbaki}, chap. 1, exercice 12
page 127). Comme $ad$ est la restriction de la représentation adjointe de l'algèbre de Lie du groupe $T_n$, les racines $\lambda$ sont à valeurs réelles.

Si $\xi\in\N$, $ad_{\xi}$ est nilpotente sur $\N$ et son image est contenue dans $\N$. Par cons{\'e}quent, la forme lin{\'e}aire $\tau:\xi\mapsto\tr ad_{\xi}$ sur $\G$ est nulle sur $\N$. Comme $G$ n'est pas unimodulaire, $\tau$ n'est pas identiquement nulle sur $\G$, donc elle est non nulle sur $\A$. 

Supposons que $\hp1(G)\not=0$. Montrons que pour tout élément $\lambda\in\Sigma$, $\lambda$, en tant que forme linéaire sur $\A$, est un multiple positif de $\tau$. Choisissons $\xi_{0}\in\A$ tel que $\tau(\xi_{0})>0$. Alors, d'apr{\`e}s le lemme \ref{strichartz}, si $\hp1(G)\not=0$, il existe une fonction non constante $w\in\mathcal{B}_{\xi_{0}}^{p}$, i.e. $w\in d^{-1}dL^{p}(G)$ est invariante par le groupe {\`a} un param{\`e}tre de translations {\`a} droite engendr{\'e} par $\xi_0$. D'après les lemmes \ref{restr} et \ref{restr'}, si $\lambda(\xi_0 )\leq \tau(\xi_0 )/p$, et si $\eta\in\N$ est un vecteur propre de $ad_{\xi_0}$ relatif à $\lambda(\xi_0 )$, alors $dw$ s'annule sur $\eta$. Le corollaire \ref{invariante'} entraîne alors que $w$ est constante, contradiction. On a donc montré que, pour tout $\xi_0 \in \A$,
\begin{eqnarray*}
\tau(\xi_0 )>0 \Rightarrow \lambda(\xi_0 )>\frac{\tau(\xi_0 )}{p}.
\end{eqnarray*}
Cela entraîne en particulier que $\tau(\xi_0 )>0 \Leftrightarrow \lambda(\xi_0 )>0$, i.e. que $\lambda$ est colinéaire à $\tau$, $\lambda=c_{\lambda}\tau$, où $\d c_{\lambda}>\frac{1}{p}$.

Pour tout $\xi\in\ker\tau$, $ad_{\xi}$ est nilpotent, donc l'idéal $\ker\tau$ est nilpotent. On peut donc remplacer la paire $(\A,\N)$ par $(\R\xi_0 ,\ker\tau)$. Autrement dit, on peut supposer que $\dim A=1$. Les formes lin{\'e}aires $\lambda\in\Sigma$ ne sont rien de plus que les nombres $\lambda=\lambda(\xi_{0})$, et leurs parties réelles sont toutes du même signe. Autrement dit, $G$ est le produit semi-direct d'un groupe de Lie nilpotent simplement connexe $N$ par $\R$ agissant par une d{\'e}rivation dont les valeurs propres ont une partie r{\'e}elle strictement positive. Un tel groupe admet une m{\'e}trique riemannienne invariante {\`a} courbure sectionnelle strictement n{\'e}gative, \cite{Heintze}.\qed 

\begin{prop}
\label{homneg}
Soit $N$ un groupe de Lie nilpotent simplement connexe de dimension $n-1$, $\alpha$ une
d{\'e}rivation de $\N$ dont les valeurs propres ont des parties r{\'e}elles
$0<\l_{1}\leq\cdots\leq\l_{n-1}$ strictement positives. Soit
$G=N\timez_{\alpha}\R$. 
Posons
$$
{\bf p}(G)=
\frac{\l_{1}+\cdots+\l_{n-1}}{\l_{1}}.
$$ 
Alors $\tp1(G)=0$ pour tout $p\geq 1$, $\rp1(G)=0$ si $p\leq {\bf p}(G)$ et $\rp1(G)\not=0$ si $p>{\bf p}(G)$.
\end{prop}

\preuve
Le r{\'e}sultat appara{\^\i}t dans \cite{P} dans le cas particulier o{\`u} $\alpha$ est semi-simple.

Le champ de vecteurs invariant {\`a} gauche $\xi_{0}$ engendrant l'action à droite du facteur $\R$ sur le produit semi-direct satisfait $ad_{\xi^{0}} =\alpha$ sur $\N$, d'où $\tr(ad_{\xi_{0}})=\l_{1}+\cdots+\l_{n-1}>0$. Le lemme \ref{strichartz} s'applique et $\hp1(G)=\mathcal{B}_{\xi_{0}}^{p}$ est s{\'e}par{\'e}. 

On a vu au paragraphe précédent que si $\hp1(G)\not=0$, alors les parties réelles $\lambda_i$ des valeurs propres de $ad_{\xi_0}$ sur $\N$ satisfont toutes $\d\l_i >\frac{\tr ad_{\xi_0}}{p}=\frac{\l_{1}+\cdots+\l_{n-1}}{p}$, soit $p>{\bf p}(G)$.

Inversement, notons $\pi:N\timez\R \to N$ la projection sur le second facteur. Si $p>{\bf p}(G)$, pour toute fonction lisse $\bar{u}$ {\`a} support compact sur $N$, $\bar{u}\circ \pi$ est dans $\mathcal{B}_{\xi_{0}}^{p}$. En effet, soit $\chi$ une fonction lisse sur $\R$ qui vaut $0$ au voisinage de $-\infty$ et $1$ au voisinage de $+\infty$. Le flot $\phi_t$ de $\xi$ donne un diff{\'e}omorphisme 
$$
N\times\R\to G,\quad (n,t)\mapsto \phi_{t}(n).
$$
Posons $u(n,t)=\chi(t)\bar{u}(n)$. Si $p>{\bf p}(G)$, $du\in L^{p}$, $u_{\infty}=\bar{u}\circ\pi\in\mathcal{B}_{\xi}^{p}$. En particulier, $\rp1(G)\not=0$.\qed

\subsection{Preuve du th{\'e}or{\`e}me \ref{degre1}}

Soit $M$ un espace homog{\`e}ne riemannien connexe non compact. Soit $G$ la composante neutre de son groupe d'i\-so\-m{\'e}\-tries. Comme $M$ et $G$ sont quasiisométriques (lemme \ref{cocompact}), $\tp1(M)=\tp1(G)$. D'apr{\`e}s le corollaire \ref{unimod}, s'il existe $p>1$ tel que $\tp1(M)\not=0$, alors $G$ est une extension compacte d'un groupe de Lie r{\'e}soluble unimodulaire. Inversement, si $G$ est une extension compacte d'un groupe de Lie r{\'e}soluble unimodulaire, alors pour tout $p>1$, $\tp1(M)\not=0$.

D{\'e}sormais, nous pouvons supposer que $\tp1(M)=0$. D'apr{\`e}s le lemme \ref{qires}, $M$ est quasiisom{\'e}trique {\`a} un groupe de Lie r{\'e}soluble $G$ qui satisfait aux hypoth{\`e}ses de la proposition \ref{nonuni}. De nouveau, $\tp1(G)=\tp1(M)=0$ donc $G$ n'est pas unimodulaire. S'il existe $p>1$ tel que $\hp1(M)\not=0$, alors $\rp1(G)\not=0$. D'apr{\`e}s la proposition \ref{nonuni}, $G$ admet une m{\'e}trique invariante {\`a} courbure sectionnelle n{\'e}gative, et le d{\'e}tail est donn{\'e} par la proposition \ref{homneg}.\qed

\subsection{Preuve du corollaire \ref{car2}}

Soit $M$ un espace homog{\`e}ne riemannien connexe non compact. L'espace des $1$-formes harmoniques $L^2$ sur $M$ est isomorphe {\`a} $R^{1,2}(M)$. J. Cheeger et M. Gromov \cite{CheegerGromov} ont montré que si $M$ est une variété riemannienne fermée à l'infini qui revêt une variété compacte, alors les nombres de Betti $L^{2}$ de $M$ (dimensions de Von Neumann des espaces de cohomologie $L^{2}$) sont tous nuls. En particulier, $R^{1,2}(M)=0$. Le même argument, fondé sur la dimension de von Neumann, s'étend sans changement au cas des groupes de Lie fermés à l'infini. Par invariance sous quasiisométrie (lemme \ref{cocompact} et corollaire \ref{h1qi}), on a donc $R^{1,2}(M)=0$ pour tout espace homogène riemannien fermé à l'infini.

D'apr{\`e}s le th{\'e}or{\`e}me \ref{degre1}, si $R^{1,2}(M)\not=0$, $M$ est commensurable {\`a} un espace homog{\`e}ne $M'$ {\`a} courbure sectionnelle strictement n{\'e}gative, et tel que ${\bf p}(M')<2$. L'expression donn{\'e}e par la proposition \ref{homneg} montre que cela n'est possible que si $\dim M'=2$.\qed

\vskip1cm
\noindent 

Univ Paris-Sud, Laboratoire de Mathématiques d'Orsay, Orsay Cedex, F-91405; \\
CNRS, Orsay cedex, F-91405\\
\smallskip
{\tt\small Pierre.Pansu@math.u-psud.fr\\
http://www.math.u-psud.fr/$\sim$pansu}

\end{document}